\documentclass[12pt]{article}
\usepackage{amssymb,amsmath}
\usepackage{hyperref}
\usepackage{graphicx}
\usepackage{color}
\usepackage{chngcntr}

\begin{document}

\title{\Large\bf Efficient computation of pi by \\
the Newton--Raphson iteration and a two-term Machin-like formula}

\author{
\normalsize\bf S. M. Abrarov\footnote{\scriptsize{Dept. Earth and Space Science and Engineering, York University, Toronto, Canada, M3J 1P3.}}
\, and B. M. Quine$^{*}$\footnote{\scriptsize{Dept. Physics and Astronomy, York University, Toronto, Canada, M3J 1P3.}}}

\date{January 3, 2018}
\maketitle

\begin{abstract}
In our recent publication we have proposed a new methodology for determination of the two-term Machin-like formula for pi with small arguments of the arctangent function of kind
$$
\frac{\pi }{4} = {2^{k - 1}}\arctan \left( {\frac{1}{{{\beta _1}}}} \right) + \arctan \left( {\frac{1}{{{\beta _2}}}} \right),
$$
where $k$ and ${\beta _1}$ are some integers and ${\beta _2}$ is a rational number, dependent upon ${\beta _1}$ and $k$. Although ${1/\left|\beta _2\right|}$ may be significantly smaller than ${1/\beta _1}$, the large numbers in the numerator and denominator of $\beta_2$ decelerate the computation. In this work we show how this problem can be effectively resolved by the Newton--Raphson iteration method.
\vspace{0.25cm}
\vspace{0.25cm}
\end{abstract}

\section{Introduction}

The two-term Machin-like formula for pi is given by
\begin{equation}\label{eq_1}
\frac{\pi }{4} = {\alpha _1}\arctan \left( {\frac{1}{{{\beta _1}}}} \right) + {\alpha _2}\arctan \left( {\frac{1}{{{\beta _2}}}} \right),
\end{equation}
where ${\alpha _1}$, ${\alpha _2}$, ${\beta _1}$ and ${\beta _2}$ are some integers or rational numbers. The Maclaurin series expansion of the arctangent function, also known as the Gregory's series \cite{Lehmer1938, Borwein1989}, can be represented as
$$
\arctan \left( x \right) = \sum\limits_{m = 1}^\infty  {\frac{{{{\left( { - 1} \right)}^{m + 1}}}}{{2m - 1}}{x^{2m - 1}}}  = x - \frac{{{x^3}}}{3} + \frac{{{x^5}}}{5} - \frac{{{x^7}}}{7} +  \ldots.
$$
Since from this series expansion it follows that $\arctan \left( x \right) = x + O\left( {{x^3}} \right)$, it is reasonable to look for a two-term Machin-like formula for pi with smaller arguments (by absolute value) of the arctangent function to improve convergence rate in computation \cite{Chien-Lih2004, Abrarov2017a, Abrarov2017b}.

Using an identity relating arctangent function with natural logarithm
$$
\arctan \left( {\frac{1}{x}} \right) = \frac{1}{{2i}}\ln \left( {\frac{{x + i}}{{x - i}}} \right),
$$
after some trivial rearrangements from equation \eqref{eq_1} it follows that \cite{Abrarov2017a}
\begin{equation}\label{eq_2}
i = {\left( {\frac{{{\beta _1} + i}}{{{\beta _1} - i}}} \right)^{{\alpha_1}}}{\left( {\frac{{{\beta _2} + i}}{{{\beta _2} - i}}} \right)^{{\alpha_2}}}.
\end{equation}

The equations \eqref{eq_1} and \eqref{eq_2} can be significantly simplified for theoretical analysis by taking ${\alpha _2} = 1$. This leads to
\begin{equation}\label{eq_3}
\frac{\pi }{4} = {\alpha _1}\arctan \left( {\frac{1}{{{\beta _1}}}} \right) + \arctan \left( {\frac{1}{{{\beta _2}}}} \right).
\end{equation}
and
$$
i = {\left( {\frac{{{\beta _1} + i}}{{{\beta _1} - i}}} \right)^{{\alpha _1}}}\frac{{{\beta _2} + i}}{{{\beta _2} - i}},
$$
respectively. The solution with respect to ${\beta _2}$ in the last equation is given by
\begin{equation}\label{eq_4}
{\beta _2} = \frac{2}{{{{\left[ {\left( {{\beta _1} + i} \right)/\left( {{\beta _1} - i} \right)} \right]}^{{\alpha _1}}} - i}} - i.
\end{equation}
It is not difficult to see that that if ${\alpha _1}$is an integer and ${\beta _1}$ is an integer or a rational number, then ${\beta _2}$ is a rational number \cite{Abrarov2017a, Abrarov2017b}.

Equation \eqref{eq_3} becomes particularly interesting by considering the second term as a small remainder defined as \cite{Abrarov2017a}
$$
\Delta  = \arctan \left( {\frac{1}{{{\beta _2}}}} \right).
$$ 
Specifically, if the condition $\left|\Delta\right| << 1$ is satisfied we can write
\begin{equation}\label{eq_5}
\frac{\pi }{4} \approx {\alpha _1}\arctan \left( {\frac{1}{{{\beta _1}}}} \right).
\end{equation}
Consequently, if there is an identity of kind
\begin{equation}\label{eq_6}
\frac{\pi }{4} = {\alpha _1}\arctan \left( {\frac{1}{\gamma }} \right),
\end{equation}
then ${\beta _1}$ in equation \eqref{eq_5} can be chosen such that ${\beta _1} \approx \gamma $ in order to obtain the term ${\alpha _1}\arctan \left( {1/{\beta _1}} \right)$ sufficiently close to $\pi /4$. This enables us to reduce the argument of the second arctangent function in equation \eqref{eq_3} such that $1/\left|\beta _2\right| << 1/\beta _1$.

It should be noted that an elegant method showing how to reduce the arguments of the arctangent function in the two-term Machin-like formula for pi was initially suggested by Chien-Lih \cite{Chien-Lih2004}. However, in contrast to the Chien-Lih's method, our approach does not require multiple step-by-step algebraic manipulations and can be developed by a relatively simple iteration procedure (see \cite{Abrarov2017b} for details).

A random choice of the values ${\alpha _1}$ and ${\beta _1}$ is inefficient. For example, substituting ${\alpha _1} = 7$ and ${\beta _1} = {10^9}$ into equation \eqref{eq_4} leads to ${\beta _2}$ that is very close to unity and, therefore, not small enough for rapid convergence (see \cite{Abrarov2017a} for specific details). However, using the following equation \cite{Abrarov2017a, Abrarov2017b}
$$
\frac{\pi }{4} = {2^{k - 1}}\arctan \left( {\frac{{\sqrt {2 - {c_{k - 1}}} }}{{{c_k}}}} \right),
$$
where the nested radicals are ${c_{k + 1}} = \sqrt {2 + {c_k}} $ and ${c_1} = \sqrt 2 $, we can construct a two-term Machin-like formula for pi providing a rapid convergence. Specifically, taking an integer ${\alpha _1} = {2^{k - 1}}$ and an integer or a rational number ${\beta _1}$ such that
$$
{\beta _1} = \frac{{{c_k}}}{{\sqrt {2 - {c_{k - 1}}} }} + \varepsilon, \qquad\qquad \left| \varepsilon  \right| << \beta_1,
$$
where $\varepsilon $ is the error term, we obtain
\begin{equation}\label{eq_7a}
\frac{\pi }{4} = {2^{k - 1}}\arctan \left( {\frac{1}{{{\beta _1}}}} \right) + \arctan \left( {\frac{2}{{{{\left[ {\left( {{\beta _1} + i} \right)/\left( {{\beta _1} - i} \right)} \right]}^{{2^{k - 1}}}} + i}} + i} \right),
\end{equation}
where
\[
\frac{2}{{{{\left[ {\left( {{\beta _1} + i} \right)/\left( {{\beta _1} - i} \right)} \right]}^{{2^{k - 1}}}} + i}} + i = \frac{1}{{{\beta _2}}},
\]
that can result in a very rapid convergence rate in computing pi (see \cite{Abrarov2017a, Abrarov2017b} for more information).

Although this implementation can provide small arguments (by absolute value) of the arctangent functions in equation \eqref{eq_3}, it leads to a value ${\beta _2}$ with a large number of the digits in its numerator and denominator that may cause some complexities in computation \cite{Guillera2017}. In this work we show how the Newton--Raphson iteration method may be applied to resolve effectively such a problem.

Several efficient algorithms for computing pi have been developed by using the Newton--Raphson iteration method (see for example \cite{Gallagher2015}). However, its application to the two-term Machin-like formula for pi may be more promising since the ratios $1/\beta_1$ and $1/\beta_2$ can be chosen arbitrarily small by absolute value in order to achieve faster convergence in computation (see \cite{Abrarov2017a} and \cite{Abrarov2017b} describing how the ratios $1/\beta_1$ and $1/\beta_2$ can be reduced by absolute value). To the best of our knowledge this approach has never been reported in scientific literature.

\section{Methodology}

There are several efficient approximations for the arctangent function \cite{Abrarov2017b, Chien-Lih2005, Milgram2005}. Recently we have derived a new series expansion of the arctangent function \cite{Abrarov2017b}
\begin{equation}\label{eq_8}
\begin{aligned}
\arctan \left( x \right) &= i\sum\limits_{m=1}^{\infty}{\frac{1}{2m-1}}\left(\frac{1}{\left(1+2i/x\right)^{2m-1}}-\frac{1}{\left(1-2i/x\right)^{2m-1}}\right) \\
&=2\sum\limits_{m = 1}^\infty  {\frac{1}{{2m - 1}}\frac{{{g_m}\left( x \right)}}{{g_m^2\left( x \right) + h_m^2\left( x \right)}}},
\end{aligned}
\end{equation}
where
$$
{g_1}\left( x \right) = 2/x, \,\, {h_1}\left( x \right) = 1,
$$
$$
{g_m}\left( x \right) = {g_{m - 1}}\left( x \right)\left( {1 - 4/{x^2}} \right) + 4{h_{m - 1}}\left( x \right)/x
$$
and
$$
{h_m}\left( x \right) = {h_{m - 1}}\left( x \right)\left( {1 - 4/{x^2}} \right) - 4{g_{m - 1}}\left( x \right)/x.
$$
Although this series expansion \eqref{eq_8} is very rapid in convergence, especially at $x <  < 1$, the large numbers in the ratio of the rational number ${\beta _2}$ decelerate the computation \cite{Guillera2017}. In order to resolve this problem we suggest the application of the Newton--Raphson iteration method.

Let
\begin{equation}\label{eq_9}
\arctan \left( x \right) = y.
\end{equation}
Then from this equation it follows that
$$
x = \tan \left( y \right)
$$
or
\begin{equation}\label{eq_10}
\tan \left( y \right) - x = 0.
\end{equation}

The Newton--Raphson iteration method is based on the formula \cite{Householder1970, Alefeld1981, Scavo1995, Recktenwald2000}

\begin{equation}\label{eq_11}	
{y_{n + 1}} = {y_n} - \frac{{f\left( {{y_n}} \right)}}{{f'\left( {{y_n}} \right)}}.
\end{equation}
Therefore, in accordance with relation \eqref{eq_10}
$$
f\left( y \right) = \tan \left( y \right) - x \Rightarrow f'\left( y \right) = \frac{d}{{dy}}\left( {\tan \left( y \right) - x} \right) = {\sec ^2}\left( y \right),
$$
the equation \eqref{eq_11} yields a very efficient iteration formula for the arctangent function
\begin{equation}\label{eq_12}
{y_{n + 1}} = {y_n} - {\cos ^2}\left( {{y_n}} \right)\left( {\tan \left( {{y_n}} \right) - x} \right),
\end{equation}
such that (see equation \eqref{eq_9})
$$
\mathop {\lim }\limits_{n \to \infty } {y_n} = \arctan \left( x \right).
$$

\section{Implementation}

The computational test shows that approximation of the function $1 - {\sin ^2}\left( {{y_n}} \right)$ is faster in convergence than an approximation of the function ${\cos ^2}\left( {{y_n}} \right)$. Therefore, taking this into account and considering that for our case $x = 1/{\beta _2}$, it is convenient to rewrite the formula \eqref{eq_12} in form
\begin{equation}\label{eq_13}
{y_{n + 1}} = {y_n} - \left( {1 - {{\sin }^2}\left( {{y_n}} \right)} \right)\left( {\tan \left( {{y_n}} \right) - \frac{1}{{{\beta _2}}}} \right).
\end{equation}

There is an expansion series for the tangent function
$$
\tan \left( \theta  \right) = \sum\limits_{m = 1}^\infty  {\frac{{{{\left( { - 1} \right)}^{m - 1}}{2^{2m}}\left( {{2^{2m}} - 1} \right){B_{2m}}{\theta ^{2m - 1}}}}{{\left( {2m} \right)!}}},	\qquad\qquad - \frac{\pi }{2} < \theta  < \frac{\pi }{2}.
$$
However, this series expansion contains the Bernoulli numbers ${B_{2m}}$ that requires intense computation when index $m$ increases. Therefore, application of this expansion series may not be optimal. Instead, we can significantly simplify the computation by rewriting equation \eqref{eq_13} in form
\begin{equation}\label{eq_14}
{y_{n + 1}} = {y_n} - \left( {1 - {{\sin }^2}\left( {{y_n}} \right)} \right)\left( {\frac{{\sin \left( {{y_n}} \right)}}{{\cos \left( {{y_n}} \right)}} - \frac{1}{{{\beta _2}}}} \right),
\end{equation}
where the sine and cosine functions can be approximated, for example, by truncating their Maclaurin series expansions as given by
$$
\sin \left( {{y_n}} \right) = \sum\limits_{m = 0}^\infty  {\frac{{{{\left( { - 1} \right)}^m}y_n^{2m + 1}}}{{\left( {2m + 1} \right)!}}}
$$
and
$$
\cos \left( {{y_n}} \right) = \sum\limits_{m = 0}^\infty  {\frac{{{{\left( { - 1} \right)}^m}y_n^{2m}}}{{\left( {2m} \right)!}}},
$$
respectively.

\section{Sample computation}

In our work \cite{Abrarov2016} we have derived a simple formula for pi
$$
\frac{\pi }{4} = {2^{k - 1}}\arctan \left( {\frac{{\sqrt {2 - {c_{k - 1}}} }}{{{c_k}}}} \right).
$$
By comparing this equation with equation \eqref{eq_6} we can see that for this specific case ${\alpha _1} = {2^{k - 1}}$ while $\gamma  = c_k/\sqrt {2 - {c_{k - 1}}}$. Consequently, we can construct the two-term Machin-like formula for pi with small arguments of the arctangent function by choosing the integer ${\beta _1}$ such that
$$
{\beta _1} \approx \frac{{{c_k}}}{{\sqrt {2 - {c_{k - 1}}} }}.
$$

Denote $\left\lfloor {\,\,\,} \right\rfloor $ as the floor function. Then the error term $\varepsilon $ can be taken as
$$
\varepsilon  = \left\lfloor {\frac{{{c_k}}}{{\sqrt {2 - {c_{k - 1}}} }}} \right\rfloor - \frac{{{c_k}}}{{\sqrt {2 - {c_{k - 1}}} }}, \qquad\Rightarrow\varepsilon < 0.
$$
Therefore, it is convenient to apply a simple equation in order to determine the integer ${\beta _1}$ as follows \cite{Abrarov2017b}
$$
{\beta _1} = \left\lfloor {\frac{{{c_k}}}{{\sqrt {2 - {c_{k - 1}}} }}} \right\rfloor.
$$

We can assign, for example, $k = 6$. Consequently, we have
$$
{\beta _1} = \left\lfloor {\frac{{\sqrt {2 + \sqrt {2 + \sqrt {2 + \sqrt {2 + \sqrt {2 + \sqrt 2 } } } } } }}{{\sqrt {2 - \sqrt {2 + \sqrt {2 + \sqrt {2 + \sqrt {2 + \sqrt 2 } } } } } }}} \right\rfloor  = 40.
$$
Substituting these values of $k$ and ${\beta _1}$ into equation \eqref{eq_7a} yields
\begin{equation}\label{eq_15}
\frac{\pi }{4} = 32\arctan \left( {\frac{1}{{40}}} \right) + \arctan \left( x \right),
\end{equation}							
where
$$
x = \frac{1}{{{\beta _2}}} =  - \frac{{{\text{38035138859000075702655846657186322249216830232319}}}}{{{\text{2634699316100146880926635665506082395762836079845121}}}}.
$$
The following is the Mathematica code that validates the Machin-like formula \eqref{eq_15} for pi by returning the output '{\bfseries{{\ttfamily{True}}}}':
\small
\begin{verbatim}
Pi/4 == 
 32*ArcTan[1/40] + 
  ArcTan[-(38035138859000075702655846657186322249216830232319/
      2634699316100146880926635665506082395762836079845121)]
\end{verbatim}
\normalsize

It is relatively easy to compute $\arctan \left( {1/\beta_1} \right)$ by using, for example, the equation \eqref{eq_8} since  ${\beta _1}$ is just an integer. Even though it is advantageous when the arguments of the arctangent function in the Machin-like formula \eqref{eq_3} for pi are smaller by absolute value, its computation may be challenging since the rational number ${\beta _2}$ consists of the numerator and denominator with large number of the digits \cite{Guillera2017}. This complexity occurs since each additional term in equation \eqref{eq_8} increases tremendously the number of the digits due to exponentiations in determination of the intermediate values ${g_m}\left( x \right)$ and ${h_m}\left( x \right)$.

Using an example based on the equation \eqref{eq_15} we can show how to overcome this problem. Suppose that the value of $\arctan \left( {1/40} \right)$ is computed with required accuracy (by equation \eqref{eq_8} or by any other approximation or numerical method) and suppose that at the beginning we know only five decimal digits of the constant pi. Then the initial value ${y_1}$ can be computed by substituting $\pi  \approx 3.14159$ into equation \eqref{eq_15}. This gives
\[
{y_1} = \frac{{3.14159}}{4} - 32\arctan \left( {\frac{1}{{40}}} \right) =  - 0.0144358958054451040550 \ldots
\]
Since we employ the approximated value of the constant pi with five decimal digits only, there is no any specific reason to compute all digits of the value ${y_1}$. Therefore, the value ${y_1}$ can be computed with same (or slightly better) accuracy than accuracy of the initial approximation $\pi  \approx 3.14159$. Thus, with only first six decimal digits
$$
{y_1} \approx  - 0.\underbrace {014435}_{6\,\,{\text{digits}}}
$$
we get
$$
\pi  \approx 4\left( {32\arctan \left( {\frac{1}{{40}}} \right) - {y_1}} \right) = 3.\underbrace{{14159}}_{5\,\,{\text{digits}}}358322178041622\ldots
$$
$5$ correct decimal digits of pi as expected. Substituting this approximated value ${y_1}$ into the Newton--Raphson iteration formula \eqref{eq_14} results in
\[
{y_2} =  - 0.014435232407997574182 \ldots
\]

Experimental observation shows that each iteration step doubles the number of the correct digits of pi \footnote{The convergence rate can be accelerated even further by using a higher order iteration like the Halley's method \cite{Alefeld1981, Scavo1995}, the Householder's method \cite{Householder1970} and so on. However, the Newton--Raphson iteration method is simplest in implementation (see for example \cite{Recktenwald2000}).}. Therefore, we can double the number of the decimal digits (from $6$ to $12$) for approximation of the value ${y_2}$ as follows
\[
{y_2} \approx  - 0.\underbrace {014435232407}_{12\,\,{\text{digits}}}.
\]
This provides
$$
\pi  \approx 4\left( {32\arctan \left( {\frac{1}{{40}}} \right) - {y_2}} \right) = 3.\underbrace {1415926535}_{10\,\,{\text{digits}}}8979323846\ldots
$$
$10$ correct decimal digits of pi. Substituting the approximated value ${y_2}$ into the Newton--Raphson iteration formula \eqref{eq_14} leads to
$$
{y_3} =  - 0.01443523240799679443929512531345 \ldots.
$$
Increasing again the number of the decimal digits by factor of two (from $12$ to $24$) for approximating $y_3$ in form
$$
{y_3} \approx  - 0.\underbrace {014435232407996794439295}_{24\,\,{\text{digits}}},
$$
we can gain
$$
\pi  \approx 4\left( {32\arctan \left( {\frac{1}{{40}}} \right) - {y_3}} \right) = 3.\underbrace {141592653589793238462643}_{24\,\,{\text{digits}}}383279 \ldots
$$
$24$ correct decimal digits of pi. This procedure can be repeated over and over again in order to achieve the required accuracy for pi.

The following is the Mathematica code showing computation of pi according to the described iteration procedure:
\bigskip
\small
\begin{verbatim}
(* M is trancating integer *)
M := 20

(* Series expansion for sine function *)
sinF[y_, M_] := Sum[((-1)^m*y^(2*m + 1))/(2*m + 1)!, {m, 0, M}]

(* Series expansion for cosine function *)
cosF[y_, M_] := Sum[((-1)^m*y^(2*m))/(2*m)!, {m, 0, M}]

(* x is argument of the arctangent function *)
x := -(38035138859000075702655846657186322249216830232319/
    2634699316100146880926635665506082395762836079845121)


Print["--------------------------------------------"]
y1 := 3.14159`200./4 - 32*ArcTan[1/40]
Print["The value of y1 is ", N[y1, 21], "..."]
y1 := -0.014435`200.
Print["Approximated value of y1 = ", N[y1, 5]]
Print["Actual value of pi is ", N[Pi, 21], "..."]
Print["Approximated value of pi is ", 
    N[4*(32*ArcTan[1/40] + y1), 21], "..."]

Print["--------------------------------------------"]
y2 := y1 - (1 - sinF[y1, M]^2)*(sinF[y1, M]/cosF[y1, M] - x)
Print["The value of y2 is ", N[y2, 20], "..."]
y2 := -0.014435232407`200.
Print["Approximated valu of y2 = ", N[y2, 11]]
Print["Actual value of pi is ", N[Pi, 21], "..."]
Print["Approximated value of pi is ", 
    N[4*(32*ArcTan[1/40] + y2), 21], "..."]

Print["--------------------------------------------"]
y3 := y2 - (1 - sinF[y2, M]^2)*(sinF[y2, M]/cosF[y2, M] - x)
Print["The value of y3 is ", N[y3, 31], "..."]
y3 := -0.014435232407996794439295`200.
Print["Approximated value of y3 = ", N[y3, 23]]
Print["Actual value of pi is ", N[Pi, 31], "..."]
Print["Approximated value of pi is ", 
    N[4*(32*ArcTan[1/40] + y3), 31]]
\end{verbatim}
\normalsize

In our publication \cite{Abrarov2017b} we have presented the two-term Machin-like formula \eqref{eq_3} for pi, where $\alpha_1 = 2^{26}$, $\beta_1 = 85445659$ and
\[
\begin{aligned}
{{\beta}_{2}}&=-\frac{\overbrace{\text{2368557598}\ldots \text{9903554561}}^{\text{522,185,816}\,\,\text{digits}}}{\underbrace{\text{9732933578}\ldots \text{4975692799}}_{\text{522,185,807}\,\,\text{digits}}} \\ 
 & =-\text{2}\text{.43354953523904089818}\ldots \times {{10}^{8}}\,\,\left( \text{rational} \right),
\end{aligned}
\]
providing $16$ decimal digits per term increment while the series expansion \eqref{eq_8} is applied \footnote{Interested reader can upload all digits of the computed rational number $\beta_2$ here: \href{https://yorkspace.library.yorku.ca/xmlui/handle/10315/33173}{https://yorkspace.library.yorku.ca/xmlui/handle/10315/33173}}. The following Mathematica code shows this convergence rate:

\bigskip
\small
\begin{verbatim}
(* Define integer k *)
k := 27

(* Define rational value beta1 *)
beta1 := 85445659

(* This is an alternative representation of equation (4), see [5] *) 
beta2 := (Cos[2^(k - 1)*ArcTan[(2*beta1)/(beta1^2 - 1)]])/(1 - 
    Sin[2^(k - 1)*ArcTan[(2*beta1)/(beta1^2 - 1)]])

(* Approximation of pi based on equations (7) and (8) *) 
piApprox[M_] := 
N[4*I*Sum[(1/(2*m - 1))*(2^(k - 1)*(1/(1 + 2*I*beta1)^(2*m - 1) - 
    1/(1 - 2*I*beta1)^(2*m - 1)) + 1/(1 + 2*I*beta2)^(2*m - 1) - 
        1/(1 - 2*I*beta2)^(2*m - 1)), {m, 1, M}], 10000] // Re

Print["Number of correct digits of pi and convergence rate:"]
Print["----------------------------------------------------"]

piDigits[M_] := Abs[MantissaExponent[Pi - piApprox[M]]][[2]]
M = 10;
While[M <= 20, {Print["At M = ", M, " number of correct digits is ", 
    piDigits[M]], Print["The convergence rate is ", 
        piDigits[M] - piDigits[M - 1], " per term increment"]}; M++]

Print["----------------------------------------------------"]
Print["Actual value of pi is"]
N[Pi, 100]

Print["At M = 5 the approximated value of pi is"]
N[piApprox[5], 100]
\end{verbatim}
\normalsize

Since the value $y_n$ in the intermediate steps of computation does not require all decimal digits, this rational number $\beta_2$ can also be approximated accordingly at each iteration. Consequently, the described application of the Newton--Raphson iteration method can be effective to resolve this problem.

\section{Conclusion}

The Newton--Raphson iteration method is applied to a two-term Machin-like formula for the high-accuracy approximation of the constant pi. The accuracy of approximated value of pi is doubled at each iteration step. Consequently, the accuracy of intermediate values $y_n$ can also be doubled stepwise at each iteration (for example by rounding or by gradual increase of the truncating integers for the relevant functions involved in computation). This methodology significantly simplifies the computation of the $\arctan \left( {1/{\beta _2}} \right)$ in the two-term Machin-like formula \eqref{eq_3} for pi and effectively resolves the problem related to the large number of the digits in the numerator and denominator of the rational value ${\beta _2}$.

\section*{Acknowledgments}

This work is supported by National Research Council Canada, Thoth Technology Inc. and York University. The authors wish to thank Dr. Jes\'us
Guillera for constructive discussions.

\bigskip

\end{document}